 \documentclass[runningheads]{llncs} 
\usepackage[T1]{fontenc}
\usepackage{amsmath,amsfonts,graphicx} 
\usepackage{tikz}
\usepackage{multirow}

\def\smallddots{\mathinner{\raise7pt\hbox{.}\raise4pt\hbox{.}\raise1pt\hbox{.}}} 
\def\smallsdots{\mathinner{\raise1pt\hbox{.}\raise4pt\hbox{.}\raise7pt\hbox{.}}}

\usepackage{algpseudocode}
\usepackage{algorithm}
\usepackage{booktabs}
\usepackage{subcaption}
\usepackage{pgfplots}
\usepackage{adjustbox}
\begin{document}
\title{New Combinations of\\ Polynomial Root-Finding Iterations}
\author{Victor Y. Pan}
\institute{
Departments of Mathematics and Computer Science \\ 
Lehman College   and \\ the Graduate Center of 
the City University of New York \\
victor.pan@lehman.cuny.edu 
 \\http://comet.lehman.cuny.edu/vpan
}
\date{} 
\maketitle   
 

\begin{abstract}
Some near-optimal polynomial root-finders of 2024-25, based on subdivision iterations, approximate all complex roots of a polynomial or all roots lying in a fixed Region of Interest in the complex plane.  We combine these iterations with Newton's  and/or Schr{\"o}der's to yield significant empirical acceleration versus each approach standing alone. 
Like the cited recent 
 algorithms, our root-finders can be applied not only to a polynomial represented  in monomial basis, by its  coefficients, but also to a {\em black box polynomial} represented by an oracle (black box subroutine) for its evaluation. 
Some by-products of our study such as     an extension of the Gauss-Lucas theorem  and a fast black box estimator
for root radius can be of independent interest.
\end{abstract}
\section{Introduction:  State of the Art and our progress}
{\bf 1.1. Polynomial root-finding: the problems.} Polynomial root-finding was studied four millennia ago by Sumerians in a rudimentary form, then stayed central in Mathematics and Computational Mathematics well into the 19th century BC 
 \cite{MB11}, and 
 is still a major subject of Computer Algebra and Symbolic-Numeric  Computing.  
 
\noindent{\bf Problem 1}: Given $\epsilon=1/2^b>0$ and real or complex coefficients  of a  polynomial 
\begin{equation}
\label{eq:p}
p=p(x): = 
p_d x^d + p_{d-1} x^{d-1}+ \cdots  + p_0 = p_d \prod_{j=1}^d (x-z_j), ~p_d \neq 0, 
\end{equation} 
 approximate within $\epsilon$ every complex {\em root} $z_j$  of $p$ and  equation $p(x)=0$. 


\noindent{\bf Problem 1$_{\mathbb S}$:} Under the assumptions of Problem 1, approximate within $\epsilon$ all roots in  a fixed convex domain (Region of Interest)
$\mathbb S$
in the complex plane.

\noindent{\bf Problems 1bb} and
{\bf  1bb$_{\mathbb S}$}
are Problems 1 and 1$_{\mathbb S}$, respectively, except that  $p$ is a black box polynomial,
represented with an oracle (black box subroutine) for its evaluation rather than its coefficients.

Black box polynomial root-finders run  faster  where a polynomial  can be evaluated  faster without involving its coefficients,  e.g.,   is a matrix polynomial,  the sum of a small number of shifted monomials 
 $(x-c)^d$,  or  a Mandelbrot-like  polynomial, defined by a recurrence (cf. \cite{B24,IP20,IP22}).
  
   One may prefer to apply
black box polynomial root-finders even to a polynomial  given with its coefficients - to  avoid the well-known heavy penalty of  coefficient
 swell for  root-squaring and for
 shifting and scaling the variable, which are basic operations for various efficient  polynomial root-finders.  
 
 Problem 1$_{\mathbb S}$ turns into Problem 1  where the domain $\mathbb S$ contains all roots. 
 
 By combining polynomial evaluation and interpolation, we can transform
Problems 1$_{\mathbb S}$ and 1bb$_{\mathbb S}$ into one another but at much higher cost than for root-finding
itself in various important applications where  the domain $\mathbb S$
contains  a small number of roots (see Appendix A). 

 \noindent{\bf 1.2. Two classes of polynomial root-finders -  brief overview.} Two classes of popular polynomial root-finders rely on (i) functional   iterations, such as Newton's  and Schr{\"o}der's  (see \cite[Ch. 5]{M07}) 
    and (ii) {\em subdivision Weyl's iterations} (see our Sec. 3);  hereafter we refer to them as 
  {\em   N, S,  and W iterations}, respectively.   
 
 N iterations converge
 with quadratic rate  to  a single simple  root locally - provided that they have been  initiated near it. \cite{BAS16,S23} formally support fast  global convergence of a variant of these iterations
 to all $d$ roots.  
 \cite{RSS24,MV24}  accelerate 
it empirically
by a factor of $d$ and  compete  with   MPSolve - user's choice library for  Problem 1.  Like 
 many other popular polynomial root-finders for Problem 1, based on functional iterations
such as  Ehrlich-Aberth's, Weirerstrass's (aka Durand -Kerner's
(see \cite[Ch. 4]{M07}), the  modified N terations  can be  applied   to Problem 1bb except that their 
 initialization  requires
 an upper bound on the root radius $r_1:=\max_{j=1}^d|z_j|$,
 whose customary  computation recipes involve
  the coefficients of $p$.  In Remark \ref{rerr}, however, we specify low cost black box computation of such a  bound and thus
  extend the  modified N terations and the other cited algorithms to Problem 1bb.
  S iterations accelerate Newton's for fast approximation of multiple roots.
     
 \cite{P24,P25,PGLZ25} apply W iterations to solve Problems 1 and 1$_\mathbb  S$ by using a near-optimal  randomized number of bit-operations and are record fast for Problems 1bb  and  1bb$_\mathbb  S$ as well. Implementations of their  early variants in 
 \cite{IP20,IP22} have competed with MPSolve 
 and superseded it for large degree polynomials.

 \noindent {\bf 1.3. Our contribution.} 
We  empirically accelerate
the bottleneck stage of W iterative root-finders of \cite{P24,P25,PGLZ25} for Problems 1bb, 1$_{\mathbb S}$,  and  1bb$_{\mathbb S}$ by means of incorporation of S   and/or  N iterations. Our  root-finders can be applied to a polynomial represented  in monomial basis as well as  to a  black box polynomial.
Some by-products of our study can be of independent interests, in particular our 
novel extension of the Gauss-Lucas theorem in Sec. \ref{sgl} and a fast black box estimator for root radius in Remark \ref{rerr}.
  
 \noindent{\bf 1.4. Organization of the paper.}  
 We devote the next section to background. 
 In Sec. \ref{ssbd} we  recall classical W iterations and their recent near-optimal acceleration.  In Sec. \ref{scmprs} we furtehr  accelerate these root-finders
 empirically - by  combining them  with S and/or N
 iterations. 
We extend the Gauss-Lucas theorem in Sec. \ref{sgl}.
In Appendix A we recall some  important applications of root-finders in the regions with small numbers of roots. In Appendix B we briefly recall the history of devising record fast root-finders. In Appendix C we prove some basic results for S
 iteration.
 
    

\section{Background}\label{sdef}
 \begin{itemize}
 \item
  $D(c,\rho):=\{x:~|x-c|\le\rho\}$ (closed disc),
 $C(c,\rho):=\{x:~|x-c|=\rho\}$ (circle),
 $A(c,\rho,\rho'):=\{x:~\rho\le |x-c|\le \rho'\}$ 
(annulus, i.e. ring).
  \item
 
 ``MCD" stands for
 ``minimal covering disc".
  A disc $D(c, \rho)$ is
$\eta$-rigid  (see Fig. 1) for $\eta>0$
if MCD of its root set
has radius at least $\eta\rho$. The maximal $\eta$ for which the disc is $\eta$-rigid is said to be its {\em rigidity}.

\begin{figure}
\centering
\resizebox{!}{0.2\textheight}{
\begin{tikzpicture}
\draw (0, 0) circle (5cm);
\draw[red] (0, 0) circle (4cm);
\draw (0, 0) circle (3cm);

\draw [->] (0, 0) edge (-4, 3);
\draw [->] (0, 0) edge (0, 4);
\draw [->] (0, 0) edge (2.4, 1.8);

\draw node at (2, 1) {$\rho$};
\draw node at (0.5, 3.5) {$\rho\theta$};
\draw node at (-3.5, 3) {$\rho\theta^2$};

\draw node at (0, -0.2) {c};
\draw node at (1, -0.8) {*};
\draw node at (2.94, 0.585) {*};
\draw node at (-0.6, 0.15) {*};
\draw node at (-2.94, -0.585) {*}; 
\draw node at (3.89, 3.95) {*};
\draw node at (4.7, -2.75) {*};
\draw node at (-4.69, 1.8) {*};
\draw node at (-4.3, 3.43) {*};
\end{tikzpicture}}
\caption  {Roots are marked by asterisks. The red circle $C(c,\rho\theta)$ 
and the  disc $D(c,\rho\theta)$ are  $\theta$-isolated.
The disc $D(c,\rho\theta)$ is $1/\theta$-rigid.
}
\end{figure}

 \item

A disc $D(c, \rho)$ or 
a circle $C(c,\rho)$  is said to be
$\theta$-{\em isolated} for $p(x)$ and $\theta>1$ if the ring (annulus) $A(c,\rho, \theta\rho))$ or $A(c,\rho/\theta,\rho\theta)$, respectively, contains no roots. The maximal $\theta$ for which the disc or circle is $\theta$-isolated is said to be its {\em isolation}, denoted
 $i(D(c,\rho))$ or $i(C(c,\rho))$,  respectively.
\item%
$\#(\mathbb S)$ denotes the number of roots lying in a domain $\mathbb S$ in the complex plane.
\item%
A $\sigma$-{\bf soft exclusion/inclusion (e/i) test}, for  $p(x)$ of (\ref{eq:p}),  
 $\sigma>1$,  and the unit disc $D(0,1)$, either outputs 1 and stops if it detects    that $r_{d}\le \sigma$, that is, that $\#(D(0,\sigma))>0$, or outputs 0 and stops\footnote{The bounds $r_{d}\le \sigma$ and $r_{d}>1$ can hold simultaneously (cf. Fig. 2), but as soon as a $\sigma$-soft e/i test verifies any of them, it stops without checking the other bound.} if it detects that $r_{d}> 1$, that is,    that $\#(D(0,1))=0$.  
\end{itemize}
\begin{figure}
\centering
\resizebox{!}{0.2\textheight}{
\begin{tikzpicture}
\draw[red] (0, 0) circle (1.414);
\draw (-1, 1) -- (1, 1) -- (1, -1) -- (-1, -1) -- (-1, 1);
\draw[red] (0, 0) circle (1); 
\draw (-0.707, 0.707) -- (0.707, 0.707) -- (0.707, -0.707) -- (-0.707, -0.707) -- (-0.707, 0.707);
\draw[red] (0, 0) circle (0.01); 
\draw node at (-1.2, 0) {*};
\draw node at (0, 1.2) {*};
\draw node at (0.9, -0.9) {*};
\end{tikzpicture}}
\caption{Both    
 exclusion and $\sqrt 2$-soft inclusion criteria hold for  e/i test applied to the internal disc.}   
 \end{figure} 

 
 We will use the following simple observation.
 \begin{lemma}\label{obscchy}
The linear map, combining shift and scaling of   the variable $x$,
  \begin{equation}\label{eqshft}
x\mapsto   
\frac{x-c}{\rho},
~p(x)\mapsto  p\Big(\frac{x-c}{\rho}\Big),~D(c,\rho)\mapsto D(0,1),~C(c,\rho)\mapsto C(0,1),
\end{equation}
 preserves the
 root sets, rigidity and isolation of the  disc $D(c,\rho)$.
\end{lemma}
The map  extends
soft e/i tests from the unit disc $D(0,1)$ to any disc.
\medskip

 Based on factorization
(\ref{eq:p}) represent
{\em Newton's inverse ratio} as follows:
\begin{equation}\label{eqratio}
{\rm NIR}(x):=\frac{p'(x)}{p(x)}=\sum_{j=1}^d\frac{1}{x-z_j}.
\end{equation}
Deduce from this equation the follwing estimates for the extremal, that is, the smallest and largest, {\em root radii}:
 \begin{equation}\label{eqrr}
r_d(x):=\min_{j=1}^d|x-z_j|\le\frac{d}{{\rm NIR}(x)}\le r_1(x):=\max_{j=1}^d|x-z_j|.
\end{equation}
Write $r_d:=r_d(0)$, $r_1:=r_1(0)$.


\section{Classical subdivision (S) root-finders and their recent acceleration}\label{ssbd}
   
\noindent{\bf 3.1. Overview.} W iterations are
traced back to  Weyl's algorithm 
of 1924,
advanced  in  \cite{H74,R87,P00} and
  used in Computational Geometry as {\em Quad-tree Construction}. They extend  bisection from root-finding in a line segment to   polynomial root-finding in a square in the complex plane, called a {\em suspect square}. 
More generally, they seek roots in a
region made up of  congruent  and pairwise disjoint, although possibly adjacent,  suspect squares that together cover
a fixed  Region  of Interest  $\mathbb S$ containing $m$ roots in its small neighborhood.   
\medskip

\noindent{\bf 3.2.  W iteration} subdivides, that is,  partitions, every suspect square into four congruent sub-squares,
 and to each 
 applies an  e/i test --
 to decide whether the sub-square contains any root lying in the Region $\mathbb S$. If  it does not, then the sub-square is discarded;
 otherwise it is called a
 suspect square and undergoes subdivision
(see Fig. 3).
\begin{figure} 
\centering
\resizebox{!}{0.2\textheight}{
\begin{tikzpicture}
\draw (-4, -4) -- (-4, 4) -- (4, 4) -- (4, -4) -- (-4, -4);
\draw (-4, 0) -- (4, 0);
\draw (0, -4) -- (0, 4);
\draw (-4, 2) -- (0, 2);
\draw (0, -2) -- (4, -2);
\draw (-2, 4) -- (-2, 0);
\draw (2, 0) -- (2, -4);
\draw (-1, 4) -- (-1, 2);
\draw (-2, 3) -- (0, 3); 
\draw (0, -3) -- (2, -3);
\draw (1, -2) -- (1, -4);
\draw node at (-1.2, 2.6) {*};
\draw node at (-1.7, 2.2) {*};
\draw node at (1.3, -2.6) {*};
\draw node at (1.7, -2.3) {*};
\end{tikzpicture}}
\caption{Four roots marked by asterisks lie in  suspect sub-squares; the empty sub-squares are discarded.}
\end{figure}

W iteration never discards a square containing a root, 
 and  at any W iteration  every root
 in $\mathbb S$ is approximated by the center of a suspect square within half-diameter, itself halved in every iteration and hence decreased by a factor of $2^b$ in $b$ 
iterations.
A root can lie in  at  most $ 4$ suspect squares. If  an initial suspect
square has diameter $2^h$, then in  at most $4(b+h)m$ e/i
tests the centers of at  most  $4m$ suspect squares approximate
within  $1/2^b$ 
all roots in the Region  $\mathbb S$.

\begin{remark}\label{reini} For every $i$, all suspect  squares 
 processed at the $i$th W iteration lie in the $\beta_i$-neighborhood of $\mathbb S$; together with  the diameter of a suspect square, the value $\beta_i$ fast decreases from $\beta_0$ to 0 as  $i$ increases to $\infty$ (see Fig. 4).
\end{remark} 

 \begin{figure}[t!]
\centering
  \begin{tabular}[c]{ccccc}
  
\begin{subfigure}[b]{0.3\textwidth}
\centering
\resizebox{!}{0.15\textheight}{
\begin{tikzpicture}
\draw[color=black!60, line width=0.1pt] (0,0) rectangle (8,8);
 \fill[red!30] (0,0) rectangle (8,8);
  \fill[white] (0,1.3)--(4, 8)--(8,1.3) --(0,1.3);
\draw[step=4, color=black!60, line width=0.1pt] (0,0) grid (8,8);
\draw[line width=1pt] (0,1.3)--(4, 8)--(8,1.3) --(0,1.3);
\end{tikzpicture}}
\subcaption{1st iteration}
\end{subfigure}

&

\begin{subfigure}[b]{0.2\textwidth}
\centering
\resizebox{!}{0.15\textheight}{
\begin{tikzpicture}
\draw[color=black!60, line width=0.1pt] (0,0) rectangle (8,8);
\fill[red!30] (0,0) rectangle (8,6);
\fill[red!30] (2,6) rectangle (6,8);
 \fill[white] (0,1.3)--(4, 8)--(8,1.3) --(0,1.3);
\draw[step=2, color=black!60, line width=0.1pt] (0,0) grid (8,8);
\draw[line width=1pt] (0,1.3)--(4, 8)--(8,1.3) --(0,1.3);
\end{tikzpicture}}
\subcaption{2nd iteration}
\end{subfigure}
\\

\begin{subfigure}[b]{0.2\textwidth}
\centering
\resizebox{!}{0.15\textheight}{
\begin{tikzpicture}
\draw[color=black!60, line width=0.1pt] (0,0) rectangle (8,8);
\fill[red!30] (0,1) rectangle (8,3);
\fill[red!30] (1,3) rectangle (7,5);
\fill[red!30] (2,5) rectangle (6,7);
\fill[red!30] (3,7) rectangle (5,8);
 \fill[white] (0,1.3)--(4, 8)--(8,1.3) --(0,1.3);
\draw[step=2, color=black!60, line width=0.1pt] (0,0) grid (8,8);
\draw[step=1, color=black!60, line width=0.1pt] (0,0) grid (8,6);
\draw[step=1, color=black!60, line width=0.1pt] (2,6) grid (6,8);
\draw[line width=1pt] (0,1.3)--(4, 8)--(8,1.3) --(0,1.3);
\end{tikzpicture}}
\subcaption{3rd iteration~~}
\end{subfigure}

& 

\begin{subfigure}[b]{0.2\textwidth}
\centering
\resizebox{!}{0.15\textheight}{
\begin{tikzpicture}
\draw[color=black!60, line width=0.1pt] (0,0) rectangle (8,8);
\fill[red!30] (0,1) rectangle (8,2.5);
\fill[red!30] (0.5,2.5) rectangle (7.5,3);
\fill[red!30] (1,3) rectangle (7,4);
\fill[red!30] (1.5,4) rectangle (6.5,5);
\fill[red!30] (2,5) rectangle (6,5.5);
\fill[red!30] (2.5,5) rectangle (5.5,6.5);
\fill[red!30] (3,6.5) rectangle (5,7.5);
\fill[red!30] (3.5,6.5) rectangle (4.5,8);
 \fill[white] (0,1.3)--(4, 8)--(8,1.3) --(0,1.3);
\draw[step=2, color=black!60, line width=0.1pt] (0,0) grid (8,8);
\draw[step=1, color=black!60, line width=0.1pt] (0,0) grid (8,6);
\draw[step=1, color=black!60, line width=0.1pt] (2,6) grid (6,8);
\draw[step=0.5, color=black!60, line width=0.1pt] (0,1) grid (8,3);
\draw[step=0.5, color=black!60, line width=0.1pt] (1,3) grid (7,5);
\draw[step=0.5, color=black!60, line width=0.1pt] (2,5) grid (6,7);
\draw[step=0.5, color=black!60, line width=0.1pt] (3,7) grid (5,8);

\draw[line width=1pt] (0,1.3)--(4, 8)--(8,1.3) --(0,1.3);
\end{tikzpicture}}
\subcaption{4th iteration} 
\end{subfigure}
\end{tabular}
\caption{$\beta_i$-neighborhoods of a triangular region shrink as $i$ increases.}
\label{fg2}
\end{figure} 
   

\noindent {\bf 3.3. Soft e/i tests and faster root-finding.}
Efficient
e/i tests are known for discs rather than squares. Moreover, an e/i test is hard if  a root  lies on or near the boundary circle, and $\sigma$-{\em soft e/i tests for  $\sigma>1$} address these problems: we claim {\em exclusion} and  discard  the
tested square if 
a soft e/i test of a disc covering this square outputs 0; we  claim $\sigma$-{\em soft inclusion} and call the square suspect otherwise.


 \begin{theorem}\label{theitst} {\rm \cite[Prop. 3.4]{P24}.}
For any real $\sigma>1$, integer $q>m$, and disc $D$  one can reduce a $\sigma$-soft e/i test of that disc to the evaluation of 
NIR$(x)$ at $2q\lceil\log_{\sigma}(1+6m\sqrt{q}))\rceil$
points.
\end{theorem}


 \begin{corollary}\label{coeitst} $O(m)$ soft e/i tests  evaluate  NIR at
$O(m^2)$ points.
\end{corollary} 

\begin{remark}\label{rerr} [Root radius via e/i tests.]
Define the {\em reverse polynomial}, 
\begin{equation}\label{eqpolyrev} 
p_{\rm rev}:=x^dp\Big(\frac{1}{x}\Big)~{\rm for}~x\neq 0;~p_{\rm rev}=p_0\prod_{j=1}^d\Big(x-\frac{1}{z_j}\Big)~{\rm for}~p_0=p(0)\neq 0,
\end{equation} 
and notice that $r_1\le 1/\rho$
if a soft e/i test (e.g., one from
\cite{P24,P25,PGLZ25})  applied to
the disc $D(0,\rho)$ and the polynomial $p_{\rm rev}$ outputs 0.  By repeating such a test recursively for fast decreasing positive 
$\rho$ until 0 is output, we  arrive at a  simple fast black box estimator for $r_1$, in contrast to the known ones involving the coefficients of $p$.  
\end{remark}

 \cite{P24} has estimated that 
 its W iterations, based on fast soft e/i tests of Thm. \ref{theitst},  approximate all roots in 
$\mathbb S$ by using $\tilde O(b^2(m^2+d))$ bit-operations 
for Problem 1$_{\mathbb S}$.  This only exceeds the 
lower bound  $0.25b\max\{d+1,(m+1)m\}$ of \cite[Cor. B.3]{P24} 
 by a factor of $\tilde O(b)$. 

\begin{figure}
\centering
\resizebox{!}{0.2\textheight}{
\begin{tikzpicture}
\draw (-2, 2) -- (2,2) -- (2, -2) -- (-2, -2) -- (-2, 2);

\draw (-2, 0) -- (2, 0);
\draw (0, 2) -- (0, -2);

\draw (-2, -1) -- (0, -1);
\draw (-1, 0) -- (-1, -2);

\draw (-2, -1.5 ) -- (-1, -1.5);
\draw (-1.5, -1) -- (-1.5, -2 );

\draw (-1.25, -1.5) -- (-1.25, -2);
\draw (-1.5, -1.75) -- (-1, -1.75);

\draw (-1.375, -1.75) -- (-1.375, -2);
\draw (-1.5, -1.875) -- (-1.25, -1.875);

\draw[red] (-1.75, 2.90) node [left]{$*$};
\draw[red] (-1.15, -1.65) node [left]{$*$};
\draw[green] (-1.23, -1.80) node [left]{$*$};      

\draw[blue] (2.225, 2) node [left]{$*$};
\draw[blue] (0.225, 0) node [left]{$*$};
\draw[blue] (-0.775, -1) node [left]{$*$};
\draw[blue] (-1.275, -1.5) node [left]{$*$};
\draw[blue] (-1.0, -1.750) node [left]{$*$};
\end{tikzpicture}}
\caption  {Five blue marked centers of suspect squares converge to a green-marked root cluster with linear rate; two red marked   
 W iterates converge with a superlinear rate.}
\end{figure}

\begin{figure}
\centering
\resizebox{!}{0.2\textheight}{
\begin{tikzpicture}
\draw (0, 0) circle (0.74);
\draw (0, 0) circle (1.42);
\draw (-2, 2) -- (2,2) -- (2, -2) -- (-2, -2) -- (-2, 2);
\draw (-0.5, 0.5) -- (0.5,0.5) -- (0.5, -0.5) -- (-0.5, -0.5) -- (-0.5, 0.5);
\draw (-1, 2) -- (-1, -2);
\draw (1, 2) -- (1, -2);
\draw (2, -1) -- (-2, -1);
\draw (2, 1) -- (-2, 1);
\draw (0, 2) -- (0, -2);
\draw (-2, 0) -- (2, 0);
\draw[red] (-0.1, 0.15) node [left]{$*$}; 
\draw[blue] (-0.1, -0.15) node [left]{$*$};
\draw[red] (-0.1, 0.15) node [right]{$*$}; 
\draw[blue] (-0.1, -0.15) node [right]{$*$}; 
\end{tikzpicture}}
\caption{A W step fast decreases the diameter of a compact component  and  of its MCD, does not decrease the  distance from the MCD to external roots, and hence  fast increases its  isolation.}
\label{fg3}
\end{figure}

\noindent {\bf 3.4.  Acceleration to near-optimal  W  iteration root-finders} has been outlined in \cite[Sec. 1.6.7]{P24} and elaborated upon in \cite{P25}. To see a {\em major obstacle} to such acceleration,  partition the union of the suspect squares processed at a W iteration into components.
Progress in root-finding
is only slow where a component
stays single and unbroken
in many consecutive W iterations 
(see Fig. 5), but then very soon the component becomes strongly isolated.
\begin{lemma}\label{leisl0}
 Any compact component made up of suspect squares
is broken into more than  one components in $O(\log(d))$ W iterations  unless (i) it is made up of  at most four suspect squares and
(ii) the iterations make it $d^{\nu}$-isolated
for a fixed constant $\nu>0$.
\end{lemma}
\begin{proof}
The union of five  suspect squares contains two roots separated by at least the side length $\lambda_0$ at the initial W iteration.
In a compact  component they must be connected with a chain of at most $d$ suspect squares
and therefore are separated by at most $d\lambda_i=d\lambda_0/2^i$ at the $i$th W iteration. 
Hence in $ \lceil\log_2(d+1)\rceil$  steps these roots must lie in distinct components. 
This proves claim (i), and one can similarly prove the next lemma. In the legend of Fig. \ref{fg3} we prove claim (ii). 
\end{proof}
\begin{lemma}\label{leisl}
   $O(\log(d))$ W iterations
applied to a minimal covering square of an $1/d^{\phi}$-rigid disc for a constant $\phi$ output more than one compact component. 
\end{lemma}
Based on these lemmas, 
 the authors of \cite{P24,PGLZ25}  proved
 correctness of their algorithms for  compression of  strongly
 isolated disc unit disc $D(0,1)$, and hence of any  strongly
 isolated disc (cf. ({\ref {eqshft})), without losing any root into  an $\eta$-rigid disc $D$  for reasonably bounded $1/\eta$. \\
 {\bf Problem C of disc compression:}
 For real constants $\alpha_0\ge 1$, $\alpha_1\ge 1$, $\beta_0\ge 0$, and $\beta_1\ge 0$ and the unit disc $D(0,1)$
containing  precisely $m$ roots  $z_1,\dots,z_m$ and  $\theta$-isolated for  $\theta=\alpha_0 d^{\beta_0}$, write $ D(c,\rho):={\rm MCD}(\{z_1,\dots,z_m\})$ and
compute a complex $y$ and a positive $\mu$ such  that 
\begin{equation}\label{eqrgt}
\mu/\rho\le \alpha_1 d^{\beta_1},
\end{equation}
\begin{equation}\label{eqincl}
\#(D(y,\mu))=m.
\end{equation}
Having solved  Problem C,  \cite{P24,PGLZ25}  apply S steps  to  a minimal covering square $S$ of the disc $D(y,\mu)$. Then $O(\log(d))$ steps must partition this single compact component into
at least two isolated compact components  
by virtue of Lemma \ref{leisl}. 
For a set of $m$ roots, the number of such partitions  and hence  the overall number of compression steps can be at most $m-1$. For polynomials with nested sets of root clusters,
this can lead to a nested process of combining steps of compression and subdivision, 
but  \cite{P25} estimated that
its overall complexity  is near-optimal.   
 
 Compression steps still tend to be the bottleneck  - they can involve order of $m\log(b)$ soft e/i tests overall, versus $O(m)$ other e/i tests in the near-optimal root-finders of \cite{P24,PGLZ25}. 
In view of Cor. \ref{coeitst}, the root-finders evaluate NIR at
order of $m^2\log(b)$ points at compression steps versus $O(m^2)$ 
points at all other steps. 

\section{Disc compression based on S and N iterations}\label{scmprs}
{\bf 4.1. Basic algorithm.}
Consider Problem C of disc compression
for $\theta=r^a$,  $a>2$, and a sufficiently large $r$
and  apply 
  S iteration 
$y:= 
x-\frac{m} {{\rm NIR}(x)}~{\rm for}~x\in C(0,r)$
for $|x|=r$. Then,
by virtue 
of Cor. \ref{coschr} of Appendix C, 
$|y-\frac{\bar s_1}{m}|\le\frac{2}{r-2}+\gamma$ for  $\bar s_1:=\frac{1}{m}\sum_{j=1}^mz_j\in D(c,\rho)$ and
$\gamma\le\frac{d-m}{m^2}\cdot\frac{(r+1)^2}{r^a-1}\cdot\frac{1}{1-\frac{d-m}{m}\cdot\frac{r+1}{r^a-1}}$. 
 
\noindent Hence $y$ lies in or near the MCD $D(c,\rho)$ for large values of $r$, and to solve Problem C it remains to compute $\mu$  satisfying  (\ref{eqrgt})
and (\ref{eqincl}).
 We apply a  low cost algorithm from \cite{P24,P25,PGLZ25} (said to be {\bf Alg. V}) to verify (\ref{eqincl}). If (\ref{eqincl}) holds for
 $\mu\le 1/2^b$,  then  $y$ approximates all $m$ roots in $D(c,\rho)$
 within $1/2^b$, and so we only need to  search for $\mu $   in the range $[1/2^b,|y|+1]$. Binary search for  the  exponent 
  $\log_2(\mu)$  in the range $[-b,\log_2(|y|+1)]$ requires  at most $\log_2(b+\log_2(|y|+1))/\log_2(E)$ invocations of Alg. V for $E$ denoting required upper bound on the relative error of approximation to $\mu$; we can let $E:=d^{\phi}/2$ for Problem C. \\
{\bf 4.2. Some lower bounds on the output error.} We begin with simple
  \begin{lemma}\label{lelbnd}
  Under the assumptions of Problem C, let $|y_1|\le 2 $ and 
  $|y_2|\le 2$
for two complex values $y_1$ and $y_2$.
Then $|y_1-y_2|-r_d(y_1)-r_d(y_2)\le 2\rho$.
\end{lemma}
  Combine the lemma with the
  bound $r_d(x)\le d/|{\rm NIR}(x)|$ 
  (shift the origin to $x$ to extend (\ref{eqrr}) for $h=1$ from $x=0$ to any $x$)  and obtain
  \begin{corollary}\label{colbnd}
 Under the assumptions of Lemma \ref{lelbnd}, it holds that
 $|y_1-y_2|-d/{\rm  NIR}(y_1)-d/{\rm  NIR}(y_2)\le 2\rho$.
  \end{corollary}
  \begin{figure}
  \label{fgNEW}
  \end{figure}
  Now we propose the following algorithm.  
  \\
{\bf Algorithm E: Estimating the radius $\rho$ of the MCD from below.}\\
INPUT: two real numbers $a>2$ and  $r>2$ such that the unit disc $D(0,1)$ contains $m$ roots $z_1,\dots,z_m$ and  is $r^a$-isolated from the other  $d-m$ roots.\\
 OUTPUT: $\bar r\le \rho$. \\
 INITIALIZATION: Fix an integer $q>0$ [e.g., $q=8$ or $q=16$]  and $q$ equally-spaced 
points $x_1,\dots,x_q$ of the circle $C(0,r)$.\\
 COMPUTATIONS:\\
1. Apply  S  iteration 
$y_i=x_i-m/{\rm NIR}(x_i)$
for $i=1,\dots,q$. \\
2. Compute  the values                                                                                                                                                                                                                                                                                                                                                                                                                                                                                                                                                                                                                                                                                                                                                                                                                                                                                                                                                                                                                                                                                                                                                                                                                                                                                                                                                                                                                                                                                                                                                                                                                                                                                                                                                                                                                                                                                                                                                                                                                                                                                                                                                                                                                                                                                                                                                                                                                                                                                                                                                                                                                                                                                                                                                                                                                                                                                                                                                                                                                                                                                                                                                                                                                                                                                                                                                                                                                                                                                                                                                                                                                                                                                                                                                                                                                                                                                                                                                                                                                                                                        
$r'(y_i):=d/|{\rm NIR}(y_i)|$ 
for  $i=1,\dots,q$. \\
3. Compute and output the value $\bar r=0.5
\max_{i,j=1}^q \{|y_i-y_j|-d/{\rm  NIR}(y_i)-d/{\rm  NIR}(y_j)\}$.                                                                                                       

{\em Correctness} of the  algorithm follows from
Cor.  \ref{colbnd}.
\\
{\bf 4.3. The first heuristic remedies.} The lower bound $\bar r$ on $1/2^b$ can  be poor  (and even negative)  if  \\(i) $\max_{i\neq j}|y_i-y_j|$ is small
and/or \\
 (ii) $d/{\rm  NIR}(y_i)\gg r_d(y_i)$ for all $i$.
 
Both 
deficiencies (i) and (ii) are severe for worst case pairs of polynomials $p$ and complex $x$.  For a heuristic remedy we increase of the number $q$ of orbits.
 
\begin{remark}\label{redstnc}
For an additional or alternative remedy against deficiency (ii) we can refine the estimate $d/$NIR$(y_i)$ for $r_{d}(y_i)$ by means of 
 modification of  stage 2  as follows:
$r_d'(y_i):=\gamma +\min_{\gamma\in\bar \Gamma}d/|{\rm  NIR}(y_i-\gamma)|$. Here $\bar \Gamma$ denotes a fixed   small finite set of absolutely small complex values, e.g., equally-spaced on the circle $C(y_i,\mu)$ for a fixed small $\mu>0$.
For another refinement of that estimate,
we can shift $y_i$ into the origin  and then compute the same bound for the polynomial $p_{2^h}(y_i)=\prod{j=1}^d(x-z_j^{2^h})$  for a fixed $h\ge 1$,
obtained from $p(x)$ in $h$ steps of 
root-squaring \cite{H59}.
In terms of the computational cost, both refinements, however,
can very well be  inferior to the  refinement by means of
the increase
of the number $q$ of orbits,
which  is a remedy for both deficiencies (i) and (ii). 
\end{remark}
\noindent{\bf 4.4. Modified
 S iteration.} Towards further increase of the value $\bar r$, we can
 modify stage 1 of Alg. E by modifying  S iteration
as follows,
\begin{equation}\label{eqemq}
y_i :=
x_i-\frac{(1-\beta)m} {{\rm NIR}(x_i)}.
\end{equation}
Here we either fix a reasonably small positive  $\beta$  based on statistics of  preliminary experiments
or define it dynamically:
e.g., initially fix $\beta=1/20$, say, and then increase it if  $\#(D(y,\bar r))<m$
or decrease it if  $\#(D(y,\bar r))=m$ and
$\bar r\ll \min_{i=1}^q|y_i-y|$.\\ 
{\bf 4.5.  Why such perturbation 
slows down S iteration.}
Our modification of Alg. E is semi-heuristic, with the goal of slowing down slightly the jump of
S  iteration into or close to the MCD $D(c,\rho)$. We do this
by means  of perturbing S iteration  with N iteration
$y= 
x-\frac{1} {{\rm NIR}(x)}$. They  also  converge to the MCD but much slower.
 To see why so,   write
  \begin{equation}\label{eqf}
 f(x)=\prod_{j=1}^m(x-z_j),~
 {\rm NIR}_f(x):=\frac{f'(x)}{f(x)},~\Delta(x):={\rm NIR}(x)-{\rm NIR}_f(x),
 \end{equation}
 recall (\ref{eqratio})
 and obtain that $\Delta(x)=\sum_{j=m+1}^d\frac{1}{x-z_j}$ and hence
  \begin{equation}\label{eqDlt}
 |\Delta(x)|\le\frac{d-m}{r^a-r},~
 |{\rm NIR}_f(x)|\ge\frac{m}{r+1} ~
 {\rm for}~|x|=r.
  \end{equation}
  Now let  $|x|=r\mapsto \infty$ to simulate N iterations in case of large $r$. Then $|\Delta(x)|\mapsto 0$, ${\rm NIR}_f(x)\mapsto
  {\rm NIR}(x)$,   the MCD $D(c,\rho)$ collapses into a single root of  multiplicity $m$, and we
  recall  that N iterations
converge 
to such a root
rather slowly - with linear rate (see \cite[Sec. 4.1]{M07}). 
\\
{\bf 4.6. Modified N iterations.}    
We have devised Alg. ${\mathcal A}$ based on   S  iteration but can replace it with  path-lifting or path following modifications of N iterations (see \cite{KS94,HSS01,BAS16,SS17,S23,RSS24,MV24}).  
Root approximations
computed  with these iterations are made up of  
  $q$ values $v_{g,h}=p(x_{g,h})$, $g=1,\dots,q$;  $h=0,1,\dots$,  
 converging to 0, while the  values $x_{g,h+1}$
 are recovered from $v_{g,h+1}=p(x_{g,h+1})$ via N type expressions. [In \cite[Thm. 5.1]{KS94} they are specified in terms of path-lifting as follows,  $x_{g,h+1}:=x_{g,h}-\frac{p(x_{g,h})-v_{g,h+1}}{p'(x_{g,h})}$.] 
 
 \cite{BAS16,S23} proved
 fast global convergence
 of these iterations to all $d$ roots, although their complexity estimates greatly exceed those of \cite{P24,P25,PGLZ25}. 
 
 Instead of root-finding task we only need to solve a much simpler Problem C of
 disc compression, and the proofs of \cite{BAS16,S23}
 are simplifed accordingly, except that we must
    extend the Gauss-Lucas theorem,
 which we do in the next section.  This extension implies that MCD $D(c,\rho)$ covers all  $m-1$ roots of $p'(x)$ (critical points)
that lie in the disc $D(0,r)$, and thus implies that convergence of orbits to these  $m-1$ attractors
would still  support disc compression.
 The extension counters the impact of $d-m$ roots lying outside the  disc  $D(0,r)$, while in the Gauss-Lucas classical case, where $m=d$,  this probem
disappears. 
 
 Presently, testing disc compression based on   modified N iterations is 
 complicated because the codes of the most efficient  variants of these iterations in \cite{SS17,RSS24,MV24}  are not publicly available.  
 
\section{Extension of Gauss-Lucas theorem}\label{scmptrcmpcmpn}\label{sgl}   
  
We begin with the following lemma.
\begin{lemma}\label{le1}
Let $|y|\le 1<z$ for a complex number $y$ and a real number $z$.
Then $\Re(\frac{1}{z-y})\ge \frac{1}{z+1}$. 
\end{lemma}

\begin{figure}[h]
\centering
\resizebox{!}{0.2\textheight}{
\begin{tikzpicture}
  \draw (-2.25,0) -- (2.25,0);
  \draw (0,-1.25) -- (0,1.25);
  
 \draw (0,0) circle [radius=1] ;
 \draw[black, fill=black] (0.7,0.7) circle (0.02); 

  \draw[above right, outer sep=-2pt, ] node at (1.0,0) {\tiny 1};
 \draw[above left, outer sep=-2pt, ] node at (0,1.0) {\tiny {\bf i}};
 \draw[above left, outer sep=-2pt, ] node at (-1.0,0) {\tiny -1};
 \draw[below left, outer sep=-2pt, ] node at (0,-1.0) {\tiny -{\bf i}};

\draw[above right, outer sep=-1pt, ] node at (0.7,0.7) {\tiny $y$};
 
  \draw[black, fill=black] (0.7,0) circle (0.02); 
 \draw[below, outer sep=-1pt] node at (0.7,0) {\tiny $\Re(y)$};
\draw (0.7,0.7) -- (0.7,0); 

 \draw[black, fill=black] (2,0) circle (0.02); 
\draw[above, outer sep=-1pt, ] node at (2,0) {\tiny $z$};

\draw (-1.0, 0) -- (0.7,0.7) -- (2,0); 

\end{tikzpicture}
}
\caption{}
\label{fig1}
\end{figure}

\begin{proof} 
Let $\phi=\phi(z-y,z)$, $0\le \phi \le \pi/2$,
denote the angle between the real axis  and the line 
passing through $z$ and $z-y$ (see Fig. \ref{fig1}) and obtain 
in polar coordinates
that
$z-y=|z-y|~(\cos(\phi)+{\bf i}~\sin(\phi))$
and furthermore, 
$\frac{1}{\cos(\phi)+{\bf i}~\sin(\phi)}=\cos(\phi)-{\bf i}~\sin(\phi)$. Hence 
$\frac{1}{z-y}=\frac{\cos(\phi)-{\bf i}~\sin(\phi)}{|z-y|},~{\rm and~ so}~
\Re\Big(\frac{1}{z-y}\Big)=\frac{\cos(\phi)}{|z-y|}.$
 Finally substitute the  bound $|z-y|\le(z+1)\cos(\phi)$ (see Fig. \ref{fig1}). 
\end{proof}



\begin{theorem}\label{thglR} 
For two integers $d$  and $m$, $1\le m\le d$, and a real $R>1$, write
$R_{1}:=
\frac{m}{d}(R+1)-1.$
For
a polynomial $p(x)$ of (\ref{eq:p})
let
\begin{equation}\label{eq1R} 
|z_j|\le 1~{\rm for}~j\le m;~|z_j|\ge R~{\rm for}~j> m.
\end{equation}
Then $m-1$ roots of the polynomial  $p'(x)$ lie in the  unit disc $D(0,1)$, while its other $d-m$  roots lie outside the open disc $\{x:~|x|<R_1\}$.
 \end{theorem}

\begin{proof}
First let  $p'(z)=0$ for $1<z<R$. 
 Hence $p(z)\neq 0$ (see (\ref{eq1R})),
and so $\frac{p'(z)}{p(z)}=0$.
Recall that
 $\frac{p'(x)}{p(x)}=\sum_{j=1}^d\frac{1}{x-z_j}$ for all complex $x$
and write
$$f(z):=\frac{p'(z)}{p(z)}=f_1(z)+ f_{2}(z),~
f_1(z):=\sum_{j=1}^m\frac{1}{z-z_j},~
 f_{2}(z):=\sum_{j=m+1}^d\frac{1}{z-z_j}.$$  (\ref{eq1R}) implies that $|f_{2}(z)|\le \frac{d-m}{R-z}$;
 Lemma \ref{le1} implies that
$|f_1(z)|\ge \Re (f_1(z))\ge \frac{m}{z+1}$. Combine these two 
bounds with the  equation $f(z)=0$
and 
obtain 
$$0=|f(z)|\ge |f_1(z)|-|f_{2}(z)|\ge \frac{m}{z+1}-\frac{d-m}{R-z}.$$
 Hence
$ 
\frac{m}{z+1}\le\frac{d-m}{R-z}$
or equivalently $(R-z)m\le(d-m)(z+1)$,
and so $z\ge R_1$ if  $ 1<z< R$.
 Extend this result to the polynomial $\bar p(x):=p(ax)$ for  $|a|=1$ replacing $p(x)$; obtain that $y\ge R_1$
if $\bar p(y)=0$ for $1<y<R$. Hence the open ring
 $\{x:~1<|x|<R_1\}$ contains no roots of $\bar p'(x)$. 

Preserve
 that property 
 in  homotopic transformation of  $p(x)$  into the polynomial $t(x):=(x+1)^m(x-R)^{d-m}$. Then the polynomials
 $p'(x)$ and
 $t'(x)=(d(x+1)-m(R+1))(x+1)^{m-1}(x-R)^{d-m-1}$ must share the numbers  $m-1$ and $d-m$ of their roots  (counted with multiplicities) that lie in the disc $D(0,1)$  and in the domain $\{x:~|x|\ge R_1\}$, respectively.
This completes the proof of Thm. \ref{thglR}.
\end{proof}
\begin{remark}\label{retght}
Can we increase the bound  $R_1$  of  the theorem? No, we cannot because  for $p(x)=t(ax)=(ax)^m(ax-R)^{d-m}$ and 
 a complex number $a$, the polynomial
 $p'(x)$
has a root $\frac{R_1}{a}$, while
the roots of $p(x)$ satisfy (\ref{eq1R}) if $|a|=1$.   
\end{remark}

\begin{corollary}\label{coglR} 
Under the assumptions of 
Thm. \ref{thglR} fix an integer $k$, $1\le k\le m$,
and  write
$$R_0:=R,~R_{h}:=
\frac{m-h}{d-h}(R_{h-1}+1)-1=\prod_{i=0}^{h-1}\frac{m-i}{d-i}(R+1)-1~{\rm for}~h=1,\dots,k.$$
Then $m-k$ roots of the polynomial  $p^{(k)}(x)$ lie in the  unit disc $D(0,1)$, while its other $d-m$  roots lie outside the open disc $\{x:~|x|<R_k\}$.
 \end{corollary}

\begin{proof}
Apply Thm. \ref{thglR} and Cor. \ref{coglR} recursively to the polynomials $p'(x)$, $p''(x)$,$\dots$ standing for $p(x)$.

 \end{proof}
\begin{remark}\label{regldsc}
Extend Thm. \ref{thglR} to any 
 disc $D(c,\rho)$
and domain $E(c,R_k\rho)$ by applying the linear map
 $$x\mapsto   
\frac{x-c}{\rho},~p(x)\mapsto  t(x)=p\Big(\frac{x-c}{\rho}\Big),~D(c,\rho)\mapsto  D(0,1),~E(c,R_k\rho)\mapsto  E(0,R_k).
$$
For $k=1$ and $R=\infty$ or $m=d$
this extension turns into the Gauss-Lucas theorem applied to a disc. 
\end{remark}




APPENDIX
\medskip

\noindent {\bf A.  Two applications.}
   \medskip

(i) Solution of  Problems 0 and 1  can be extended to root-finding for a triangular multivariate  polynomial  system of equations  \cite{CLOS97,IPY21}.
The overall number of solutions of such a system can grow exponentially with its total degree, while only a small number of them lie in the Region of User's  Interest.\footnote{In computer aided design a user might want to design a 3D object under some fixed geometric constraints translated into a  system of polynomials equations. This system can have an enormous number of solutions, mostly
spurious, and the user can greatly benefit from restriction of  the computations
 to  a fixed Region of Interest.}
(ii) {\em Real root-finding} (in  a  
 real line segment)  is highly important because in applications, e.g., to optimization in algebraic geometry and geometric modeling,  
   only real roots  of a polynomial are of interest, and  typically they are much less numerous than all $d$ complex roots  \cite{K43,EK95}.
   \medskip
   
   {\bf B. Brief history and related works.} 
      \medskip
      
 Hundreds of efficient root-finders have appeared \cite{M07,MP13} and keep appearing.
In 1982, in  \cite{S82}, Sch{\"o}nhage proposed to compare their efficiency by estimating  bit-operations involved, reaching $\tilde O(bd^3$) in his  record fast  solution of Problem 1.
  Subsequent extensive research   has culminated in  
  1995 (see \cite{P95,P02}) with  bound 
$\tilde O((b+d))d^2)$  -  near-optimal for $d=\tilde O(b)$.  The next  near-optimal algorithms for  Problems 1  and 1$_{\mathbb  S}$ also solved Problems 1bb and 1bb$_{\mathbb  S}$ (see \cite{P24,P25,PGLZ25}) record fast.  \cite{IP20,IP22} implemented and tested their initial variants from \cite{Pa}. Based on distinct novel techniques, the
advanced papers \cite{BSSY18,M21,IM23} deduced similar cost bounds  for
 Problem 1 explicitly or implicitly simplified with  imposing lower bounds $\delta>0$ on the pairwise distances between roots of $p$ (compare a  discussion in \cite[Secs. 19 and 20]{S82}). This was a major restriction because the input and rounding errors turn multiple roots  into tiny root clusters. The estimates for the bit operation complexity of 
 the efficient root-finders of  \cite{BSSY18} and  \cite{M21,IM23} 
 include the terms 
$d\log(1/\delta)$ and  
$d\log(\kappa(\delta))$, 
respectively, where
$\kappa(\delta)$ denotes
the maximal condition  of the roots;  both of these terms grow to
$\infty$ as $\delta$ decreases to 0.\footnote{The announcement  in
\cite{BSSXY16} of  a new near-optimal solution of Problem 1 (with no restriction on $\delta$) stirred confusion because  for all proofs  \cite{BSSXY16}
 referred to \cite{BSSY18}.
The  best implementation -  in \cite{IPY18} - of the algorithm of \cite{BSSXY16}  was 
greatly superseded already by the implementations in \cite{IP20,IP22} of some root-finders of \cite{Pa}.}

These and all other known polynomial root-finders preceding \cite{Pa,P24,P25,PGLZ25}  involve coefficients of $p$.  The only exception is the advanced pioneering root-finder of
\cite{LV16}, although it  only approximates a single absolutely largest root of a  polynomial
 having only real roots.

{\bf C. Computing  a point in MCD by means of
S iteration.}

\begin{theorem}\label{thschr}
  For an integer $r>\max_{j=1}^d|z_j|$ and complex $x$, such that  $|x|>2r$, let
\begin{equation}\label{eqschr}
y:=x-\frac{d}{{\rm NIR}(x)}. 
\end{equation}
  Then 
$\Big|y-\frac{s_1}{d}\Big|\le \frac{2r^{2}}
{|x|-2r}=\frac{2r}{(|x|/r)-2}.$
\end{theorem}
\begin{proof}
First let $r=1$.
Combine
(\ref{eqratio})
and (\ref{eqschr}) and deduce that \\ 
$y=x-\frac{d}{{\rm NIR}(x)}=x-\frac{d}{\sum_{j=1}^d\frac{1}{x-z_j}}=x-
\frac{xd}{\sum_{j=1}^d\frac{1}{1-z_j/x}}.$\\
Recall that
$|\frac{z_j}{x}|<1$ for all $j$, by theorem's assumptions for $r=1$, substitute Newman's expansion $\frac{1}{1-v}=\sum_{g=0}^{\infty}v^g$ for 
$v=\frac{z_j}{x}$ into above equation, and obtain  $$\frac{y}{x}=1-\frac{d}{\sum_{j=1}^d\sum_{g=0}^{\infty}z_j^g/x^g}=1-\frac{d}{\sum_{g=0}^{\infty}\sum_{j=1}^dz_j^g/x^g}.$$
Substitute equations
$s_g=\sum_{j=1}^dz_j^g$,  $g=1,2,\dots$, and deduce that
$$\frac{y}{x}=1-\frac{d}{d+\sum_{g=1}^{\infty}s_g/x^g}=1-\frac{1}{1+u}~{\rm for}~u=\sum_{g=1}^{\infty}\frac{s_g}{dx^g}.$$
Notice that 
$1-\frac{1}{1+u}=
u-\frac{u^2}{1+u}$, and so 
$y=xu-\frac{xu^2}{1+u}.$
Furthermore,
$
xu=\frac{s_1}{d}+\Gamma~{\rm for}~ \Gamma:=\sum_{g=2}^{\infty}\frac{s_g}{d}x^{1-g}.$ 
Therefore, 
\begin{equation}\label{eqx}
y-\frac{s_1}{d}=\Gamma-\frac{xu^2}{1+u},~{\rm and~so}~\Big|y-\frac{s_1}{d}\Big|\le |\Gamma|+\Big|\frac{xu^2}{1+u}\Big|.
\end{equation}
Notice that
$|s_h|/d\le 1$ for all $h$
because $|z_j|\le 1$ for all $j$ and
 deduce that
 $|u|\le \frac{1}{|x|}\sum_{g=0}^{\infty}\frac{1}{|x|^g}\le \frac{1}{|x|-1}$
 and hence
$$|xu^2|\le \frac{|x|}{(|x|-1)^2};~\frac{1}{|1+u|}\le\frac{1}{1-\frac{1}{|x|-1}}=\frac{|x|-1}{|x|-2}.$$
 Furthermore,
  $$|\Gamma|\le\sum_{g=2}^{\infty}|x|^{1-g}=
\frac{1}{|x|-1}.$$
Combine these bounds with Eqn. (\ref{eqx}) and obtain Thm.  \ref{thschr} in the case of $r=1$:
$$\Big|y-\frac{s_1}{d}\Big|\le \frac{1}{|x|-1}+\frac{|x|}{(|x|-1)(|x|-2)}=\frac{2}{|x|-2}.$$  

Next scale the variable $x\mapsto x/r$ implying that $D(0,r) \mapsto D(0,1)$,
$y \mapsto y/r$ (cf. (\ref{eqratio})), and $s_1 \mapsto s_1/r$,
extend Thm.  \ref{thschr}  to the case of any $r>0$.
\end{proof}
\begin{corollary}\label{coschr}
For real values  $a>2$ and $r>2$, let
S iteration
be applied at a point $x\in C(0,r)$ towards compression of $r^a$-isolated unit disc $D(0,1)$ containing precisely 
roots $z_1,\dots,z_m$, for $1<m\le d$. Then
the iteration outputs a complex value $y$ 
such that
$$\Big|y-\frac{\bar s_1}{m}\Big|\le\frac{2}{r-2}+\gamma~{\rm for}~\gamma\le\frac{d-m}{m^2}\cdot\frac{(r+1)^2}{r^a-1}\cdot\frac{1}{1-\frac{d-m}{m}\cdot\frac{r+1}{r^a-1}}.$$
\end{corollary}
\begin{proof}
Apply Thm. \ref{thschr}
to the polynomial
$f(x)=\prod_{j=1}^m(x-z_j)$ and deduce that  $|\bar y-\frac{\bar s_1}{m}|\le \frac{2}{r-2}$ for $\bar y:=x-m/$NIR$_f(x)$, NIR$_f(x)=\frac{1}{\sum_{j=1}^m(x-z_j)}$, $y-\bar y=\frac{m}{{\rm NIR}(x)}-\frac{m}{{\rm NIR}_f(x)}=\frac{m\Delta(x)}{{\rm NIR}(x)\cdot{\rm NIR}_f(x)}$ for $\Delta(x)={\rm NIR}(x)-{\rm NIR}_f(x)$ of (\ref{eqf}).

Recall that $|\Delta(x)|\le \frac{d-m}{r^a-1}$ (see (\ref{eqDlt})),
notice that $|{\rm NIR}_f(x)|\ge \frac{m}{r+1}$,
and $|{\rm NIR}(x)|\ge |{\rm NIR}_f(x)|-|\Delta(x)|\ge \frac{m}{r+1}-\frac{d-m}{r^a-1}$. Hence
$|y-\bar y|\le\frac{d-m}{r^a-1}\cdot \frac{r+1}{m}\cdot\psi$ for $\psi:=\frac{1}{\frac{m}{r+1}-\frac{d-m}{r^a-1}}=\frac{r+1}{m}\cdot\frac{1}{1-\frac{d-m}{m}\cdot\frac{r+1}{r^a-1}}$,  and  we arrive at the corollary.
\end{proof}
\begin{remark}\label{reschr}
We can strengthen the bounds of Thm. \ref{thschr} and Cor. \ref{coschr} if we apply
S iteration to the polynomial $p_{2^h}(x)$ of Remark \ref{rerr}. One should choose a small positive $h$ to keep computational cost reasonable. 
\end{remark}

\end{document}